\documentclass[final]{siamltex}

\usepackage[parfill]{parskip}    % Activate to begin paragraphs with an empty line rather than an indent
\usepackage{amsmath,amssymb,latexsym,enumerate,mathrsfs}
\usepackage{hyperref}
\usepackage{array}
\usepackage{times}
\usepackage[]{mdframed}
\usepackage{epstopdf}
\usepackage{subfig}
\usepackage{graphicx}
\usepackage{hyperref} 
\usepackage{multicol}
\usepackage{framed}
\usepackage{moreverb}
\usepackage{stackengine}
\usepackage{marvosym}
\usepackage{color,soul}
\definecolor{lightblue}{rgb}{.90,.95,1}
\sethlcolor{yellow}

\DeclareMathOperator*{\argmax}{arg\,max}
\newtheorem{remark}{Remark}

\title{On characterizing optimal learning trajectories in a class of learning problems}

\author{Getachew K. Befekadu}

\begin{document}
\maketitle

\renewcommand{\thefootnote}{\arabic{footnote}}

\begin{abstract}
In this brief paper, we provide a mathematical framework that exploits the relationship between the maximum principle and dynamic programming for characterizing optimal learning trajectories in a class of learning problem, which is related to point estimations for modeling of high-dimensional nonlinear functions. Here, such characterization for the optimal learning trajectories is associated with the solution of an optimal control problem for a weakly-controlled gradient system with small parameters, whose time-evolution is guided by a model training dataset and its perturbed version, while the optimization problem consists of a cost functional that summarizes how to gauge the quality/performance of the estimated model parameters at a certain fixed final time w.r.t. a model validating dataset. Moreover, using a successive Galerkin approximation method, we provide an algorithmic recipe how to construct the corresponding optimal learning trajectories leading to the optimal estimated model parameters for such a class of learning problem.
\end{abstract}


\section{Statement of the problem} \label{S1}
In this section, we formally present our problem statement and describe the underlying core concepts that will allow us to provide a framework based on the relationship between the maximum principle and dynamic programming for characterizing the optimal learning trajectories in a class of learning problem.\footnote{Some numerical works using the proposed frameworks have been done and detailed results will be presented elsewhere.}

Here, we consider a class of learning problem, in which the learning algorithm's generalization performance is associated with an optimal control problem
\begin{align}
J^{\epsilon}[u] &= \Phi\bigl(\theta^{\epsilon}(T), \mathcal{Z}^{(2)}\bigr) \quad \to \quad \min_{u(t) \in U ~ \text{a.e. in } ~ t \in [0, T]} \label{Eq1.1}
\end{align}
over all solutions of the following weakly-controlled gradient system with small parameters
\begin{align}
 \dot{\theta}^{\epsilon}(t) &= - \nabla J_{0} \bigl(\theta^{\epsilon}(t), \mathcal{Z}^{(1)}\bigr) + \epsilon D\bigl(\theta^{\epsilon}(t), \tilde{\mathcal{Z}}^{(1)}\bigr) u(t), ~ \text{a.e. in } ~ t \in [0, T], ~ \theta^{\epsilon}(0)=\theta_0, \label{Eq1.2}
\end{align}
where the problem statement consists of the following core concepts and general assumptions:
\begin{enumerate} [(a).]
\item {\bf Datasets}: We are given two datasets, i.e., $\mathcal{Z}^{(k)} = \bigl\{ (x_i^{(k)}, y_i^{(k)})\bigr\}_{i=1}^{m_k}$, each with data size of $m_k$, for $k=1, 2$. These datasets, i.e., $\mathcal{Z}^{(1)}$ and $\mathcal{Z}^{(2)}$, may be generated from a given original dataset $\mathcal{Z}^{(0)} =\bigl\{ (x_i^{(0)}, y_i^{(0)})\bigr\}_{i=1}^{m_0}$ by means of bootstrapping with/without replacement. Here, we assume that the first dataset $\mathcal{Z}^{(1)}=\bigl\{ (x_i^{(1)}, y_i^{(1)})\bigr\}_{i=1}^{m_1}$ will be used for model training purpose, while the second dataset $\mathcal{Z}^{(2)} = \bigl\{ (x_i^{(2)}, y_i^{(2)})\bigr\}_{i=1}^{m_2}$ will be used for evaluating the quality of the estimated model parameter. Moreover, the perturbed dataset $\tilde{\mathcal{Z}}^{(1)}$ (which is associated with the nonlinear term $D$ in the system dynamics) is obtained by adding small random noise, i.e., $\tilde{\mathcal{Z}}^{(1)} = \bigl\{ (x_i^{(1)}, \tilde{y}_i^{(1)})\bigr\}_{i=1}^{m_1}$, with $\tilde{y}_i^{(1)} = y_i^{(1)} + \varepsilon_i$ and $\varepsilon_i \sim \mathcal{N}(0, \sigma^2)$, for $i=1,2, \ldots, m_1$, with small variance $\sigma^2$.
\item {\bf Learning algorithm via weakly-controlled gradient systems with small parameters}: We are tasked to find for a parameter $\theta \in \Theta$, from a finite-dimensional parameter space $\mathbb{R}^p$ (i.e., $\Theta \subset \mathbb{R}^p$), such that the function $h_{\theta}(x) \in \mathcal{H}$, i.e., from a given class of hypothesis function space $\mathcal{H}$, describes best the corresponding model training dataset as well as predicts well with reasonable expectation on a different model validating dataset. Here, the search for an optimal parameter $\theta^{\ast} \in \Theta \subset \mathbb{R}^p$ can be associated with a weakly controlled-gradient system of Equation~\eqref{Eq1.2}, whose {\it time-evolution} is guided by the model training dataset $\mathcal{Z}^{(1)}$ and its perturbed version $\tilde{\mathcal{Z}}^{(1)}$, i.e.,
\begin{align*}
 \dot{\theta}^{\epsilon}(t) = - \nabla J_{0} \bigl(\theta^{\epsilon}(t), \mathcal{Z}^{(1)}\bigr) + \epsilon  D\bigl(\theta^{\epsilon}(t), \tilde{\mathcal{Z}}^{(1)}\bigr)u(t), 
 \end{align*}
 where $J_0\bigl(\theta, \mathcal{Z}^{(1)}\bigr) = \frac{1}{m_1} \sum\nolimits_{i=1}^{m_1} {\ell} \bigl(h_{\theta}(x_i^{(1)}), y_i^{(1)} \bigr)$, and $\ell$ is a suitable loss function that quantifies the lack-of-fit between the model and the datasets. Moreover, $u(t)$ is a real-valued admissible control function from a compact set $U \subset \mathbb{R}^p$ entering into the system dynamics through the nonlinear term $D$.\footnote{Note that the control function $u(t)$ is admissible if it is measurable and $u(t) \in U \subset \mathbb{R}^p$ for all $t \in [0, T]$.} The parameter $\epsilon$ is a small positive number and the nonlinear term $D$ is given by
\begin{align*}
 &D\bigl(\theta, \tilde{\mathcal{Z}}^{(1)}\bigr) =\\
  & \quad \quad \operatorname{diag}\bigl \{\bigl(\partial J_0 \bigl(\theta, \tilde{\mathcal{Z}}^{(1)}\bigr)/ \partial \theta_1\bigr)^2, \bigl(\partial J_0 \bigl(\theta, \tilde{\mathcal{Z}}^{(1)}\bigr)/ \partial\theta_2\bigr)^2,\ldots, \bigl(\partial J_0 \bigl(\theta, \tilde{\mathcal{Z}}^{(1)}\bigr)/ \partial \theta_p\bigr)^2 \bigr\}.
\end{align*}
Note that the small random noise $\varepsilon_i \sim \mathcal{N}(0, \sigma^2)$, for $i=1,2, \ldots, m_1$, with small variance $\sigma^2$, in the dataset $\tilde{\mathcal{Z}}^{(1)}$ will provide a {\it dithering effect}\footnote{For example, we can set the standard deviation $\sigma = c \Vert y_i \Vert_{\infty}$ of the maximal training dataset, where $c > 0$ is a weighted noise level such as $1\%$, $5\%$ or $10\%$.}, i.e., causing some distortion to the model training dataset $\mathcal{Z}^{(1)}$ so that the control $u(t)$ will have more effect on the learning dynamics.\footnote{In Equation~\eqref{Eq1.2} above, the weakly controlled-gradient system with small parameter can be expressed as follows:
\begin{align*}
 \dot{\theta}_i^{\epsilon}(t) = - \partial J_{0} \bigl(\theta^{\epsilon}(t), \mathcal{Z}^{(1)}\bigr)/ \partial \theta_i + \epsilon \bigl(\partial J_0 \bigl(\theta, \tilde{\mathcal{Z}}^{(1)}\bigr)/ \partial \theta_i\bigr)^2 u_i(t), \quad i = 1,2, \ldots, p,
\end{align*}
where the $i^{\rm th}$-control $u_i(t)$, with $i = 1,2, \ldots, p$, enters into the system dynamics through the nonlinear term $\epsilon \bigl(\partial J_0 \bigl(\theta, \tilde{\mathcal{Z}}^{(1)}\bigr)/ \partial \theta_i\bigr)^2$ (see \cite{r1} for additional discussions, but in the context of perturbation theory).}
\item {\bf Optimal control problem}: For a given $\epsilon \in (0, \epsilon_{\rm max})$, determine an admissible optimal control $u(t) \in U$, a.e. in $t \in [0, T]$, that minimizes the following cost functional
 \begin{align*}
 J^{\epsilon}[u] = \Phi\bigl(\theta^{\epsilon}(T), \mathcal{Z}^{(2)}\bigr), \quad \text{s.t. ~~~ Equation}~\eqref{Eq1.2},
\end{align*}
 where $\Phi\bigl(\theta, \mathcal{Z}^{(2)}\bigr)$ is a scalar function that depends on the model validating dataset $\mathcal{Z}^{(2)}$. Note that the optimal control problem together with the above weakly-controlled gradient system provides a mathematical apparatus how to improve generalization performance in such a class of learning problems.\footnote{In this paper, we consider the following locally Lipschitz cost function
 \begin{align*}
\Phi\bigl(\theta, \mathcal{Z}^{(2)}\bigr) = \big(1/m_2\big) \sum\nolimits_{i=1}^{m_2} {\ell} \bigl(h_{\theta}(x_i^{(2)}), y_i^{(2)} \bigr), ~~ \text{w.r.t. the model training dataset} ~ \mathcal{Z}^{(2)},
\end{align*}
as a cost functional $J^{\epsilon}[u]$ that serves as a measure for evaluating the quality of the estimated optimal parameter $\theta^{\ast} = \theta^{\epsilon}(T)$, i.e., when $\theta^{\epsilon}(t)$ is evaluated at a certain fixed time $T$.}
 \item {\bf General assumptions:} Throughout this paper, we assume the following conditions: (i) the set $U$ is compact in $\mathbb{R}^p$ and the final time $T$ is fixed, (ii) the function $\Phi\bigl(\theta, \mathcal{Z}^{(2)}\bigr)$ is locally Lipschitz, and (iii) for any $\epsilon \in (0, \epsilon_{\rm max})$ and all admissible bounded controls $u(t)$ from $U$ a.e. $t \in [0, T]$, the solution ${\theta}^{\epsilon}(t)$ for all $t \in [0, T]$ which corresponds to the weakly-controlled gradient system of Equation~\eqref{Eq1.2} starting from an initial condition $\theta^{\epsilon}(0)=\theta_0$, exists and bounded.\footnote{For such an optimal control problem, these assumptions are sufficient for the existence of a nonempty compact reachable set $\mathcal{R}(\theta_0) \subset \Theta$, for some admissible controls on $[0,T]$ that belongs to $U$, starting from an initial point $\theta^{\epsilon}(0)=\theta_0$ (e.g., see \cite{r2} for related discussions on the Filippov's theorem providing a sufficient condition for compactness of the reachable set).}.
 \end{enumerate}

In what follows, we assume that there exists an admissible optimal control $u(t) \in U$ a.e. in $t \in [0, T]$ and for $\epsilon \in (0, \epsilon_{\rm max})$. Then, the necessary optimality conditions for the optimal control problem with weakly-controlled gradient system satisfy the following Euler-Lagrange critical point equations
\begin{align}
 \dot{\theta}^{\epsilon}(t) &= \frac{\partial H^{\epsilon}\bigl(\theta^{\epsilon}(t), p^{\epsilon}(t), u(t)\bigr)}{\partial p}, \notag \\
                                       &= - \nabla J_{0} \bigl(\theta^{\epsilon}(t), \mathcal{Z}^{(1)}\bigr) + \epsilon  D\bigl(\theta^{\epsilon}(t), \tilde{\mathcal{Z}}^{(1)}\bigr) u(t),  \quad \theta^{\epsilon}(0)=\theta_0, \label{Eq1.3} \\
   \dot{p}^{\epsilon}(t) &= -\frac{\partial H^{\epsilon}\bigl(\theta^{\epsilon}(t), p^{\epsilon}(t), u(t)\bigr)}{\partial \theta}, \notag \\
                                       &= \nabla^2 J_{0} \bigl(\theta^{\epsilon}(t), \mathcal{Z}^{(1)}\bigr) p^{\epsilon}(t) - \epsilon   u^T(t) \nabla^2 J_{0}\bigl(\theta^{\epsilon}(t), \tilde{\mathcal{Z}}^{(1)} \bigr) D \bigl(\theta^{\epsilon}(t), \tilde{\mathcal{Z}}^{(1)}\bigr) p^{\epsilon}(t), \notag\\
                                       & \quad \quad \quad \quad \quad  \quad p^{\epsilon}(T) = - \nabla \Phi\bigl(\theta^{\epsilon}(T), \mathcal{Z}^{(2)}\bigr), \label{Eq1.4}\\
                                       u(t) &= \argmax H^{\epsilon}\bigl(\theta^{\epsilon}(t), p^{\epsilon}(t), u(t)\bigr), ~~ u(t) \in U ~~ \text{a.e. in} ~~ t \in [0, T], \label{Eq1.5}
\end{align}
where the Hamiltonian function $H^{\epsilon}$ is given by
\begin{align}
 H^{\epsilon} \bigl(\theta, p, u \bigr) = \bigl \langle p,\, - \nabla J_{0} \bigl(\theta, \mathcal{Z}^{(1)}\bigr) + \epsilon  D\bigl(\theta, \tilde{\mathcal{Z}}^{(1)}\bigr)u \bigr \rangle \label{Eq1.6}
\end{align}
and such optimality conditions are the direct consequence of the {\it Pontryagin's maximum principle} (e.g., see \cite{r3} for additional discussions on the first-order necessary optimality conditions; see also \cite{r4} for related discussions in the context of learning). In the following section, we attempt to relate the above Euler-Lagrange critical point equations with that of the dynamic programming principle for characterizing the optimal learning trajectories for such a class of learning problem. Moreover, using a Galerkin method, we provide approximation solutions corresponding the optimal learning trajectories leading to the optimal estimated model parameters.

\section{Main results} \label{S2}
In this section, we exploit the relationship between the maximum principle and that of the dynamic programming principle for characterizing the optimal learning trajectories in a class of learning problem. First, let us restate the above Euler-Lagrange critical point equations to have a more useful form that will allow us to characterize the optimal trajectories for the optimization problem in Equation~\eqref{Eq1.1} with a weakly-controlled gradient system of Equation~\eqref{Eq1.2} leading to the optimal estimated model parameters.

\begin{proposition}
Suppose that the general assumptions in Section~\ref{S1} hold true. Then, the trajectory pair $\bigl(\bar{\theta}^{\epsilon}, \bar{u} \bigr)$ of the weakly-controlled gradient system with a small parameter in Equation~\eqref{Eq1.2}, with $\bar{\theta}^{\epsilon}(0)=\theta_0$, is optimal for the optimization problem in Equation~\eqref{Eq1.1} if and only if the solution for the adjoint system equation $p^{\epsilon} \colon [0, T] \to \mathbb{R}^p$, i.e.,
 \begin{align}
    -\dot{p}^{\epsilon}(t) &= \bigl \langle -\nabla^2 J_{0} \bigl(\bar{\theta}^{\epsilon}(t), \mathcal{Z}^{(1)}\bigr) + \epsilon \bar{u}^T(t) \nabla^2 J_{0}\bigl(\bar{\theta}^{\epsilon}(t), \tilde{\mathcal{Z}^{(1)}} \bigr) D \bigl(\bar{\theta}^{\epsilon}(t), \tilde{\mathcal{Z}}^{(1)}\bigr), \,p^{\epsilon}(t) \bigr \rangle, \notag\\
                                       & \quad \quad \quad \quad \quad  \quad p^{\epsilon}(T) = - \nabla \Phi\bigl(\bar{\theta}^{\epsilon}(T), \mathcal{Z}^{(2)}\bigr)\label{Eq2.1}
\end{align}
satisfies the maximum principle
\begin{align}
   & \bigl(p^{\epsilon}(t),\, - \nabla J_{0} \bigl(\bar{\theta}^{\epsilon}(t), \mathcal{Z}^{(1)}\bigr) + \epsilon D\bigl(\bar{\theta}^{\epsilon}(t), \tilde{\mathcal{Z}}^{(1)}\bigr) \bar{u}(t) \bigr)= \notag\\
     &\quad \sup_{u(t) \in U} \biggl \langle -\nabla^2 J_{0} \bigl(\bar{\theta}^{\epsilon}(t), \mathcal{Z}^{(1)}\bigr) + \epsilon  D \bigl(\bar{\theta}^{\epsilon}(t), \tilde{\mathcal{Z}}^{(1)}\bigr) \nabla^2 J_{0}\bigl(\bar{\theta}^{\epsilon}(t), \tilde{\mathcal{Z}^{(1)}} \bigr) u(t), \,p^{\epsilon}(t) \biggr \rangle, \notag \\
     &\quad\quad\quad\quad\quad \quad \quad \quad \quad \quad \quad \text{a.e. in} ~~ t \in [0, T]. \label{Eq2.2}
\end{align}
Moreover, the optimal estimated parameters $\theta^{\ast}$ can be recovered from
\begin{align}
 \theta^{\ast} = \bar{\theta}^{\epsilon}(T). \label{Eq2.3}
\end{align}
\end{proposition}

Note that our interest is to relate the adjoint system equation in Equation~\eqref{Eq2.1}, which appears in the maximum principle, to that of the value function $V(t, \theta)$ arising from the dynamic programming principle, which is associated with perturbations in initial time and state of optimal control problem in Equations~\eqref{Eq1.1} and \eqref{Eq1.2}. Recall that the dynamic programming principle is mainly concerned with the properties of the value function $V(t, \theta)$ and its characterization as a solution of the Hamilton-Jacobi-Bellman (HJB) partial differential equation. That is,
\begin{align}
 -p^{\epsilon}(t) = \nabla_{\theta} V(t,\bar{\theta}^{\epsilon}(t)), ~~ \text{for all} ~ t \in [0, T], \label{Eq2.4}
\end{align}
where
\begin{align}
 V(0, \bar{\theta}^{\epsilon}(0)) = \inf \bigl\{ \Phi\bigl(\bar{\theta}^{\epsilon}(T)\bigr)\, \bigl\vert \, \bar{\theta}^{\epsilon}(t) ~ \text{is the solution of Equation}~\eqref{Eq1.2} \,\, \text{on} \,\, [0, T], \,\, \bar{\theta}^{\epsilon}(0) = \theta_0 \bigr\}. \label{Eq2.5}
\end{align}
\begin{remark}
Note that Equation~\eqref{Eq2.4} holds true if the value function $V(\cdot,\,\cdot)$ is smooth and continuously differentiable. However, in most cases such an assumption may not hold (see \cite{r5} and \cite{r6} for additional discussions). On the other hand, if the value function is Lipschitz continuous for each $t \in [0, T]$, then we can replace by the differential inclusion $-p^{\epsilon}(t) \in \nabla_{\theta} V(t,\bar{\theta}^{\epsilon}(t))$ for all $t \in [0, T]$.
\end{remark}
\subsection*{A successive Galerkin approximation method} In this subsection, we provide an algorithmic recipe how to construct the corresponding optimal learning trajectories leading to the optimal estimated model parameters for such a class of learning problem. Here, our focus is to approximate the Euler-Lagrange critical point equations using a Galerkin method, so that we can solve numerically by the method of successive approximations. In order to this, we choose a set of $N$ basis functions $\left \{ \psi_j(t) \right\}_{j=1}^N$ and we assume that the admissible control $u(t)$ for all $t \in [0, T]$, can be expressed as
\begin{align*}
 u_i(t) = \sum\nolimits_{j=1}^N c_{ij} \psi_j(t), \quad i=1,2, \dots, p,
\end{align*}
where the coefficients $c_{ij} \in \mathbb{R}$ for $i=1,2, \dots, p$ and  $j=1,2, \dots, N$, are unknown and to be determined by the method of successive approximations. Note that we can rewrite the above equation in a more compact form as $u(t) = C \Psi(t)$, where $C=\bigl(c_{ij}\bigr)$, $i=1,2, \dots, p$, $j=1,2, \dots, N$, and $\Psi(t) = \bigl[\psi_1(t),\, \psi_2(t), \dots, \psi_N(t)\bigr]^T$, for $t \in [0, T]$. Suppose that $\bar{\theta}^{\epsilon}(t)$,  $p^{\epsilon}(t)$ for all $t \in [0, T]$, with an admissible optimal control $\bar{u}(t) =  \bar{C} \Psi(t)$, are associated with the maximum principle of Equation~\eqref{Eq2.2}, where $\bar{C}$ satisfies
\begin{align*}
 &\bigl(p^{\epsilon}(t),\, - \nabla J_{0} \bigl(\bar{\theta}^{\epsilon}(t), \mathcal{Z}^{(1)}\bigr) + \epsilon D\bigl(\bar{\theta}^{\epsilon}(t), \tilde{\mathcal{Z}}^{(1)}\bigr) \bar{u}(t) \bigr)= \\
     &\quad  \max_{\underset{j=1,2, \dots, N} {c_{ij} \in \mathbb{R}, ~ i=1,2, \dots, p}} \biggl \langle -\nabla^2 J_{0} \bigl(\bar{\theta}^{\epsilon}(t), \mathcal{Z}^{(1)}\bigr) + \epsilon  D \bigl(\bar{\theta}^{\epsilon}(t), \tilde{\mathcal{Z}}^{(1)}\bigr) \nabla^2 J_{0}\bigl(\bar{\theta}^{\epsilon}(t), \tilde{\mathcal{Z}^{(1)}} \bigr) C \Psi(t), \,p^{\epsilon}(t) \biggr \rangle.
\end{align*}
Finally, the following generic algorithm can be used to construct the corresponding optimal learning trajectories leading to the optimal estimated model parameters for such a class of learning problem.

{\rm \small

{\bf ALGORITHM: A successive Galerkin approximation method}
\begin{itemize}
\item[{\bf 0.}] {\bf Initialize:} $c_{ij}^0$, for all $i=1,2, \dots, p$ and $j=1,2, \dots, N$.
\item[{\bf 1.}] Using the admissible control $\bar{u}^{k}(t) = C^k \Psi(t)$, i.e., $\bar{u}_i^k(t) = \sum\nolimits_{j=1}^N c_{ij}^k \psi_i(t)$, with $i =1,2, \ldots, p$, for $t \in [0, T]$, solve both the forward and backward-equations in Equations~\eqref{Eq1.3} and \eqref{Eq1.4}.
 \item[{\bf 2.}] Then, update the admissible control $\bar{u}^{k+1}(t) = C^{k+1} \Psi(t)$ using
\begin{align*}
 \bar{u}_i^{k+1} (t) =  \sum_{j=1}^N c_{ij}^{k+1} \psi_i(t), \quad i =1,2, \ldots, p,
\end{align*} 
where 
\begin{align*}
 c_{ij}^{k+1} = c_{ij}^{k} + \gamma_{ij} \frac{\delta H^{\epsilon}} {\delta c_{ij}}, \quad j =1,2, \ldots, N
 \end{align*} 
and $\gamma_{ij} \in [0, 1]$ (see also Equations~\eqref{Eq1.5} and \eqref{Eq2.2}).
\item[{\bf 3.}] With the updated admissible control $\bar{u}^{k+1}(t)$, repeat Steps $1$ and $2$, until convergence, i.e., $\bigl\Vert \sum_{i=1}^p \sum_{j=1}^N \delta H^{\epsilon}/ \delta c_{ij} \bigr\Vert \le \epsilon_{\rm tol}$, for some error tolerance $\epsilon_{\rm tol} > 0$.\footnote{Note that the algorithm terminates, when $\Vert \bar{u}^{k+1} (t) - \bar{u}^{k} (t)\Vert < \epsilon_{\rm tol}$ which also implies $\bar{u}^{k+1}(t) \approx \bar{u}(t)$.}
\item[{\bf 10.}] {\bf Output:} Return the optimal estimated parameter value $\theta^{\ast} = \bar{\theta}^{\epsilon}(T)$.
 \end{itemize}}

\end{document}